# Modification of Pickands' Dependence Function for the Simulate Ordered Bivariate Extreme Data

### Mohd Bakri Adam*

*Mathematics Department,*
*Faculty of Science,*
*Universiti Putra Malaysia,*
*43400 UPM Serdang,*
*Selangor, MALAYSIA.*
*e-mail:* bakri@math.upm.edu.my

**Abstract:** We study the characteristics of the Pickands' dependence function for bivariate extreme distribution for minima, BEVM, when considering the stochastics ordering of the two variables. The existing Pickand's dependence function terminologies and theories are modified to suit the dependence functions of extreme cases. We successful implement and apply these functions to our simulate extreme data

## Contents

## 1. Introduction

In this article, the concentration of the study is on the bivariate extreme distribution for minima, BEVM, with a specific dependence assumption which lead to the development of restricted Pickands' dependence function. We develop bivariate extremal models and associated statistical procedures for vector observations whose components are subject to an order relationship. The joint structure of the extreme value models for ordered variables are needed to extract the information from extreme value data. We restrict our scope of study by considering the order of the two variables involved i.e. one variable is greater than the second variable.

---

*I would like to thanks the Government of Malaysia for the Scholarship under SLAB.

We elaborate the theoretical properties of the bivariate extreme value distribution for minima, BEVM with constraint of the variable $X$ is less than variable $Y$, $X < Y$. We include the dependence case with $X$ and $Y$ in the form of exponential margin, but with Fréchet margin and the generalized extreme value for minima (GEVM) margin also being used in certain places.

The main difference of this article compared to Nadarajah, (6), is in our development of the BEVM model, we are dealing with the *restricted logistic* function and GEVM marginal distribution. We introduce some *new* terminologies and theories related to the modification of dependence functions, $A$. We use Pickands non-parametric estimator of the dependence function to guide us in the analysis.

The structure of this paper is as follows. In Sections 2 and 3 we explain the theoretical background of the bivariate extreme value distribution and ordering properties of BEVM. Later in Sections 4 and 5, we discuss the asymmetry logistic function leading to a new formulation of a restricted logistic dependence function and a restricted exponential measure function. We also illustrated some features of the dependence function in Section 6. Later in the paper, we also illustrated how to estimating the boundary constant value $c$ using non-parametric estimation of dependence function. At the end of this paper, an application using simulation data is also demonstrated.

## 2. Bivariate Extreme Value

Bivariate extreme value theory can also be used to examine the dependence between two series of extreme data, $X$ and $Y$. Let $(X, Y)$ follow the bivariate extreme value distribution which arises as the limiting distribution of componentwise maxima see (3). He studies the BEV for maxima but as here we are interested in minima we will consider the BEVM distribution with exponential margins initially. The general bivariate extreme survival distribution function for random variables $X$ and $Y$ with exponential margins, is defined as

$$G(x,y) = \Pr(X > x, Y > y) = \exp[-V(x,y)] = \exp\left\{-(x+y)A\left(\frac{x}{x+y}\right)\right\}$$

where the exponential measure function

$$V(x,y) = \int_0^1 \max\left[\omega x, (1-\omega)y\right] dH(\omega)$$

with

$$\int_0^1 \omega dH(\omega) = 1 = \int_0^1 (1-\omega) dH(\omega) \tag{2.1}$$

where the characteristic distribution function $H(q)$, is a finite positive measure over $[0,1]$, and dependence function

$$A(\omega) = \int_0^1 \max\{(1-\omega)q, \omega(1-q)\}dH(q)$$
$$= \omega \int_0^\omega (1-q)dH(q) + (1-\omega)\int_\omega^1 qdH(q) \qquad (2.2)$$

where $(1-\omega)q > \omega(1-q)$ if and only if $\frac{q}{1-q} > \frac{\omega}{1-\omega}$, and hence if and only if $0 < \omega < q < 1$. Pickands, (7), defined $A(\omega)$ as the dependence function of $(X,Y)$. From equations (2.1) and (2.2) we have $A(0) = A(1) = 1$, $\max(\omega, 1-\omega) \leq A(\omega) \leq 1$ for $0 \leq \omega \leq 1$ and $A(\omega)$ is a convex function in the region $0 \leq \omega \leq 1$. The estimation of the Bivariate Extreme Value (BEV) distribution reduces to estimation of the $A$ function when we assume that the marginal distributions have already been transformed to exponential, see (10) and (9) for more features of $A(\omega)$.

The joint survivor function $G(x, y)$ is uniquely determined by specification of any of $V$, $A$ and $H$ functions. The stronger dependence between $X$ and $Y$ arises when $H$ places mass closer to the centre of $[0,1]$ in which case $A(\omega)$ is closer to the bounding curve $\max(\omega, 1-\omega)$ for $\omega \in [0,1]$. The $X$ and $Y$ are totally independent if $A(\omega) = 1$. The $X$ and $Y$ are perfectly dependent if $A(\omega) = \max(\omega, 1-\omega)$. From equation (2.2) the first derivative of $A(\omega)$ is $A'(\omega) = H(\omega) - 1$, as $\int_0^1 qdH(q) = 1$. We get the relation

$$H(\omega) = A'(\omega) + 1 \qquad (2.3)$$

and clearly $h(\omega) = H'(\omega) = A''(\omega)$.

Tawn, (10–12) explained in detail about the bivariate extreme distribution and more generally the multivariate extreme distribution. Here, we use only the logistic model because it is widely used due to it's simple structure. Some details about BEV distribution are also available in (1), (2) and (9).

A widely used non-parametric estimate of $A$, and here $V(x, y)$ and $G(x, y)$, is the Pickands' estimator see (7; 8). Pickands' estimator is defined as follow

$$\widehat{A}_p(\omega) = \frac{n}{\sum_{i=1}^n T_{\omega,i}} \qquad (2.4)$$

with $T_{\omega,i} = \min\left(\frac{X_i}{\omega}, \frac{Y_i}{(1-\omega)}\right)$. But $\widehat{A}_p(0) = \frac{1}{\bar{y}} \neq 1$ and $\widehat{A}_p(1) = \frac{1}{\bar{x}} \neq 1$. We use a modification of estimator (2.4) proposed by (4; 5) where they suggest replacing $T_{\omega,i}$ by

$$T_{\omega,i} = \min\left(\frac{X_i}{\omega \bar{X}}, \frac{Y_i}{(1-\omega)\bar{Y}}\right). \qquad (2.5)$$

This estimator satisfies $\widehat{A}_p(0) = \frac{\bar{y}}{\bar{y}} = 1$ and $\widehat{A}_p(1) = \frac{\bar{x}}{\bar{x}} = 1$ and $\widehat{A}_p(\omega) \geq \max(\omega, 1-\omega)$.

We want to develop a related statistical methods for vector observations whose components are subject to an order restriction, $X < Y$. We use an asymmetric logistic function from (6) and (10) as our dependence function. When the order restriction is accounted for we obtain a restricted logistic function, a modification from the asymmetric logistic function. We obtained the exponential measure function for the restricted case and later we study the features of the restricted dependence function.

## 3. Ordering Properties of BEVM

We get the BEVM$(\mu_x, \sigma_x, \xi_x, \mu_y, \sigma_y, \xi_y)$ distribution function from the survival function of the extreme minima data where

$$G(x,y) = P(X > x, Y > y)$$

$$= \exp\left\{-\int_0^{C(x,y)} \frac{1-\omega}{\left[1-\xi_y(\frac{y-\mu_y}{\sigma_y})\right]_+^{1/\xi_y}} dH(\omega) \right.$$

$$\left. - \int_{C(x,y)}^1 \frac{\omega}{\left[1-\xi_x(\frac{x-\mu_x}{\sigma_x})\right]_+^{1/\xi_x}} dH(\omega)\right\}$$

with $C(x,y) = \frac{1}{1+D(x,y)}$, where $D(x,y) = \frac{\left[1-\xi_y\left(\frac{y-\mu_y}{\sigma_y}\right)\right]_+^{1/\xi_y}}{\left[1-\xi_x\left(\frac{x-\mu_x}{\sigma_x}\right)\right]_+^{1/\xi_x}}$. For $x$ in the range of $Y$, $(x \in R_Y)$, when we take $Y = x$ then for $X < Y$,

$$\Pr(X > x | Y = x) = \frac{-\frac{\partial G(x,y)}{\partial y}|_{y=x}}{g(x)} = 0,$$

where $g(x) = -\frac{\partial G(\infty,x)}{\partial x}$, marginal density distribution for $Y$. As $g(x) \neq 0$, then for $\Pr(X > x | Y = x) = 0$ we need $\frac{\partial G(x,y)}{\partial x}|_{y=x} = 0$. However,

$$\frac{\partial G(x,y)}{\partial y} = \frac{-G(x,y)}{\sigma_y\left[1-\xi_x\left(\frac{y-\mu_y}{\sigma_y}\right)\right]_+^{1+1/\xi_y}} \int_0^{C(x,y)} (1-\omega) dH(\omega),$$

as

$$\frac{1-C(x,y)}{\left[1-\xi_x\left(\frac{y-\mu_y}{\sigma_y}\right)\right]_+^{1/\xi_y}} = \frac{C(x,y)}{\left[1-\xi_x\left(\frac{y-\mu_y}{\sigma_y}\right)\right]_+^{1/\xi_y}}.$$

When $y = x$,

$$\left.\frac{\partial G(x,y)}{\partial y}\right|_{y=x} = \frac{-G(x,x)}{\sigma_y\left[1-\xi_x\left(\frac{x-\mu_y}{\sigma_y}\right)\right]_+^{1+1/\xi_y}} \int_0^{C(x,x)} (1-\omega) dH(\omega).$$

So as $\Pr(X > x | Y = x) = 0$ for all $x \in R_Y$ holds only when $\int_0^{C(x,x)} (1-\omega) dH(\omega) = 0$. This implies that

$$h(\omega) = 0 \quad \text{for all } \omega \in [0, C(x,x)], \tag{3.1}$$

for a given $x \in R_Y$. As we need this to be zero for all $x$ so $h(\omega) = 0$ for all $\omega \in [0, c]$ where $c = \max_{x \in R_Y} C(x,x)$ or $c = \frac{1}{1+d}$ where $d = \min_{x \in R_Y} D(x,x)$.

## 4. Asymmetric Logistic Function

We want to develop a bivariate extreme value model and related statistical methods for vector observations whose components are subject to an order restriction, $X < Y$. We use an asymmetric logistic function from (6) and (10) as our dependence function. When the order restriction is accounted for we obtain a restricted logistic function, a modification from the asymmetric logistic function. We obtained the exponential measure function for the restricted case and later we study the features of the restricted dependence function.

Nadarajah, (6), and Tawn, (10), defined the exponential measure function for asymmetric logistic model with exponential margins as follows

$$V_a(x,y) = \int_0^1 \max[\omega x, (1-\omega)y] \, dH_a(\omega)$$

with

$$h_a(\omega) = (s-1)(\theta_1 \theta_2)^s [\omega(1-\omega)]^{s-2} [(\theta_1 \omega)^s + \theta_2^s (1-\omega)^s]^{\frac{1}{s}-2}, \tag{4.1}$$

and $H_a(\{0\}) = 1 - \theta_1$ and $H_a(\{1\}) = 1 - \theta_2$, with

$$\int_0^1 \omega dH_a(\omega) d\omega = \int_0^1 (1-\omega) dH_a(\omega) d\omega = 1.$$

But,

$$V_a(x,y) = (x+y) A_a\left(\frac{x}{x+y}\right).$$

From (10) if $h(\omega)$ is from the asymmetric logistic function then we get

$$A_a(\omega) = (1-\theta_1)\omega + (1-\theta_2)(1-\omega) + \{(\theta_1 \omega)^s + [\theta_2(1-\omega)]^s\}^{\frac{1}{s}}$$

$$A_a'(\omega) = (1-\theta_1) - (1-\theta_2) + [\theta_1^s \omega^{s-1} - \theta_2^s(1-\omega)^{s-1}][\theta_1^s \omega^s + \theta_2^s(1-\omega)^s]^{\frac{1}{s}-1}.$$

From equation (2.3) where $H(\omega) = A'(\omega) + 1$ we get

$$H_a(\omega) = 1 - \theta_1 + \theta_2 + [\theta_1^s \omega^{s-1} - \theta_2^s(1-\omega)^s][\theta_1^s \omega^s + \theta_2^s(1-\omega)^s]^{\frac{1}{s}-1}, 0 < \omega < 1$$

and then

$$H_a(0_+) = (2 - \theta_1) - (1 - \theta_2) + (-\theta_2^s)(\theta_2^s)^{\frac{1}{s}-1} = 1 - \theta_1,$$

and

$$H_a(1_-) = 2 - \theta_1 - (1 - \theta_2) + \theta_1^s \theta_1^{1-s} = 1 + \theta_2.$$

## 5. Getting a Modification Dependence Function

We use the notation $N$ for the logistic function used by (6) which is defined on $(a, b)$, with $0 \leq a \leq \frac{1}{2}, \frac{1}{2} \leq b \leq 1$, and $RL$ for our restricted function which is defined on $(c, 1)$, with $0 \leq c \leq \frac{1}{2}$.

Using Theorem 5.5 from (6) gives

$$h_N(\omega) = (s-1)(b-a)(\alpha\beta)^s \left[(\omega - c)(b - \omega)\right]^{s-2}$$
$$\times \left[\alpha^s(\omega - a)^s + \beta^s(b - \omega)^s\right]^{\frac{1}{s}-2}, \omega \in (a, b), s > 1 \quad (5.1)$$

as a valid density of the measure $H_N$ on $[a, b]$ with $H_N(\{a\}) = \gamma_1$ and $H_N(\{b\}) = \gamma_2$ and satisfying equation (2.3). Given $a < \frac{1}{2}$, $b \geq \frac{1}{2}$, $\gamma_1, \gamma_2 \geq 0$ are such that $\gamma_1 + \gamma_2 < 2$, $\gamma_2 \leq \frac{1-2a}{b-a}$ and $\gamma_1 \leq \frac{2b-1}{b-a}$, and $\alpha = 2b - 1$ and $\beta = 1 - 2a$.

From equation (3.1), $h(\omega) = 0$ for $\omega \in [0, c]$ then the solution for $V(x, y)$ is

$$V_N(x, y) = \int_0^1 \max\left[\omega x, (1-\omega)y\right] h_N(\omega) d\omega = \int_c^1 \max\left[\omega x, (1-\omega)y\right] h_N(\omega) d\omega. \quad (5.2)$$

Using equation (5.1) with $a = c$, $b = 1$, $\alpha = 1$, $\beta = 1 - 2c$, $\omega \in (c, 1)$ we focus on the restricted function so

$$h_{RL}(\omega) = (s-1)(1-c)(1-2c)^s[(\omega - c)(1 - \omega)]^{s-2}$$
$$\times [(\omega - c)^s + (1 - 2c)^s(1 - \omega)^s]^{\frac{1}{s}-2}. \quad (5.3)$$

Comparing equations (4.1), (5.1) and (5.3), we get $\theta_1 = \alpha = 1$ and $\theta_2 = \beta = 1 - 2c$ with $\gamma_1 = \gamma_2 = 0$.

We transform the function with $w = \frac{\omega - c}{1-c} \Rightarrow \omega = c + (1-c)w$, $0 < w < 1$, and $\frac{d\omega}{dw} = 1 - c$, we get then the exponential measure for the restricted logistic function as follows

$$V_{RL}(x, y) = \int_0^1 \max\left\{[c + (1-c)w]x, [(1-c)(1-w)]y\right\} \frac{h_a(w)}{1-c} dw$$
$$= \int_0^1 \max\left[\left(\frac{c}{1-c} + w\right)x, (1-w)y\right] h_a(w) dw$$

with $h_a(w)$ defined in equation (4.1), with $\theta_1 = 1$ and $\theta_2 = 1 - 2c$. It is easy to show that

$$V_{RL}(x, y) = \begin{cases} x & \text{for } 0 < \frac{y}{x+y} \leq c \\ \frac{1}{1-c}\left(\{[(1-c)y - cx]^s + (1-2c)^s x^s\}^{\frac{1}{s}} + cx\right) & \text{for } c < \frac{y}{x+y} < 1. \end{cases} \quad (5.4)$$

The proof is provided at the end of this paper.

Now we identify the characteristics of new dependence function, $A_{RL}$, of the restricted logistic dependence function by plotting $A_{RL}(w)$ versus $w$.

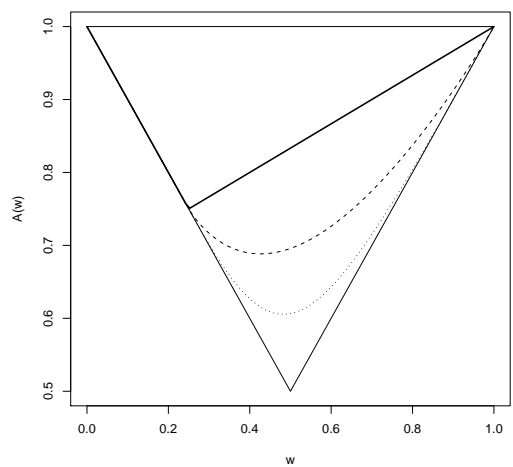

FIG 1. *We fix $c = 0.25$ but with varieties of $s = 1$ (solid line), $1.5$ (dashed line) and $2.5$ (dotted line) values.*

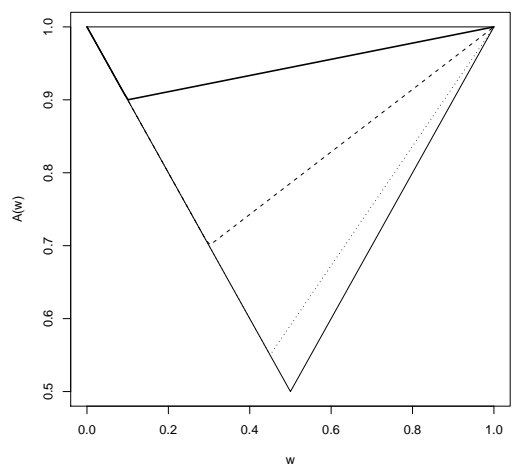

FIG 2. *We fix $s = 1$ but with varieties of $c = 0.1$ (solid line), $c = 0.3$ (dashed line) and $c = 0.45$ (dotted line) values.*

In Figure 1 we fixed $c = 0.25$ but with three values of $s$, where $s = 1.00, 1.50$ and $2.50$. When $X < Y$ with $c = 0.25$, the $X$ and $Y$ are as weakly dependent

as possible if

$$A_{RL}(w) = \begin{cases} 1 - w & \text{for} \quad 0 \leq w \leq 0.25 \\ \frac{1}{3}w + \frac{2}{3} & \text{for} \quad 0.25 < w \leq 1. \end{cases}$$

By increasing the values of $s$ to 1.50 and 2.50 we get a stronger dependence between $X$ and $Y$ in which $A(w)$ is closer to the bounding curve $\max(w, 1-w)$ for $w \in (0.25, 1]$.

In Figure 2 we study the features when $s = 1$, i.e. the weakest possible dependence between $X$ and $Y$ given the value of $c$. By increasing the values of $c$, $A(w)$ values are given in Figure 2.

## 6. Extension Cases for Ordering

In this section, we will derive the dependence function $A_N$ for the Nadarajah model when $h_N(\omega) = 0, \omega \in (0, c)$ then we propose the dependence functions, $A_N$ for when $h_N$ is only non-zero in $[0, c]$ then in $[c_1, c_2]$.

**Proposition 6.1.** *Let $A_N(\omega)$ be the dependence function of $(X, Y)$ when the measure $h_N$ is Nadarajah logistic function with $h_N(\omega) = 0$ for $\omega \in (0, c)$, $0 \leq c < \frac{1}{2}$, then $A_N(\omega)$ function is*

$$A_N(\omega) = \begin{cases} 1 - \omega & \text{if } 0 \leq \omega \leq c, \\ \frac{c(1-\omega)}{1-c} + \frac{1}{1-c}[(\omega - c)^s + (1 - 2c)^s(1 - \omega)^s]^{\frac{1}{s}} & \text{if } c < \omega \leq 1. \end{cases}$$

*with $1 \leq s$.*

**Proof**

For $0 \leq \omega \leq c$, we get

$$A_N'(\omega) = -1, \quad A_N''(\omega) = 0 = h_N(\omega).$$

For $c \leq \omega \leq 1$, we get $A_N''(\omega)$ as follow

$$(s-1)(1-c)(1-2c)^s[(\omega-c)(1-\omega)]^{s-2}[(\omega-c)^s + (1-2c)^s(1-\omega)^s]^{\frac{1}{s}-2}.$$

Compared with the density measure of equation (5.3), then $A_N''(\omega) = h_N(\omega)$.

When $X < Y$, $h(\omega) = 0$ for $0 \leq \omega \leq c$, where $0 < c < \frac{1}{2}$. But this is not the only stochastic ordering of interest in general, other orderings might be studied. Possible ordering as $X < Y^*$ where $Y^*$ can be any function such as $1/Y$, $\log Y$, $Y^2$ etc. We might also be interested in $X > Y$. We propose a possible result for new dependence functions for any ordering of $X < Y^*$ where the dependence is of a logistic type such that the following two cases follow

1. If $h_N(\omega) = 0$ for $c < \omega < 1$, where $\frac{1}{2} < c < 1$.
2. If $h_N(\omega) = 0$ for $0 < \omega < c_1$ and $c_2 < \omega < 1$ where $0 < c_1 < \frac{1}{2}$ and $\frac{1}{2} < c_2 < 1$.

For any possible relationship between $X$ and $Y$, that leads to $h_N(\omega) = 0$ for $c < \omega < 1$, where $\frac{1}{2} < c < 1$ we get

$$V_N(x,y) = \int_0^c \max(\omega x, (1-\omega)y)\, h_N(\omega)d\omega.$$

**Proposition 6.2.** *Let $A_N(\omega)$ be a dependence function of $(X,Y)$ when the measure $h_N$ is Nadarajah logistic dependence function with $h_N(\omega) = 0$ for $c < \omega < 1$, with $\frac{1}{2} < c < 1$ then $A_N(w)$ function can be written as*

$$A_N(\omega) = \begin{cases} \frac{(1-c)\omega}{c} + \frac{1}{c}[(2c-1)^s \omega^s + (c-\omega)^s]^{\frac{1}{s}} & \text{if } 0 \leq \omega \leq c, \\ \omega & \text{if } c < \omega \leq 1, \end{cases}$$

*with $1 \leq s$.*

When $h_N(\omega) = 0$ for $0 < \omega < c_1$ and $c_2 < \omega < 1$ where $0 < c_1 < \frac{1}{2}$ and $\frac{1}{2} < c_2 < 1$.

$$V_N(x,y) = \int_{c_1}^{c_2} \max(\omega x, (1-\omega)y)\, h_N(\omega)d\omega.$$

**Proposition 6.3.** *Let $A_N(\omega)$ be a dependence function of $(X,Y)$ when the measure $h_N$ is Nadarajah logistic dependence function with $h_N(\omega) = 0$ for $0 < \omega < c_1$ and $c_2 < \omega < 1$ with $0 < c_1 < \frac{1}{2} < c_2 < 1$, then $A_N(\omega)$ function is*

$$A_N(\omega) = \begin{cases} 1 - \omega & \text{if } 0 \leq \omega < a, \\ \frac{(1-c_2)(\omega-c_1)}{c_2-c_1} + \frac{c_1(c_2-\omega)}{c_1-c_2} + W(\omega) & \text{if } c_1 \leq \omega \leq c_2, \\ \omega & \text{if } c_2 < \omega \leq 1, \end{cases}$$

*where $W(\omega) = \frac{1}{c_2-c_1}[(2c_2-1)^s(\omega-c_1)^s + (1-2c_1)^s(c_2-\omega)^s]^{\frac{1}{s}}$ with $1 \leq s$.*

The proof of the latter two propositions similiar to Proposition 6.1. We show the dependence function by using the Proposition 6.3. In Figure 3 we fixed $c_1 = 0.25$ and $c_2 = 0.75$ but with three values of $s$, where $s = 1$, 2 and 5. When $s = 1$, we get the weakest possible dependence for given $c_1, c_2$ when

$$A_N(\omega) = \begin{cases} 1 - \omega & \text{if } 0 \leq \omega < c_1, \\ 0.75 & \text{if } c_1 \leq \omega \leq c_2, \\ \omega & \text{if } c_2 < \omega \leq 1. \end{cases}$$

By increasing the values of $s$ from 1 to 2 and 5 we get a stronger dependence between $X$ and $Y$ in which $A_N(\omega)$ is closer to the bounding curve $\max(\omega, 1-\omega)$ for $\omega \in (0.25, 0.75)$. In Figure 4 we show $A_N(\omega)$ for fixed $c_1 = 0.10$ and $c_2 = 0.60$, $s = 1$, 2 and 5.

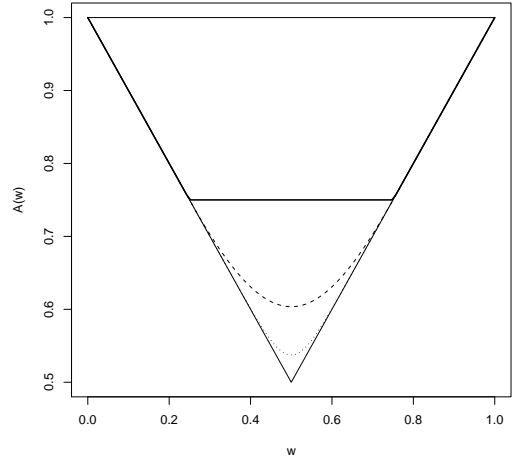

FIG 3. *We fix $a = 0.25$ and $b = 0.75$ but with varieties of $s = 1$ (solid line), $s = 2$ (dashed line) and $s = 5$ (dotted line) values.*

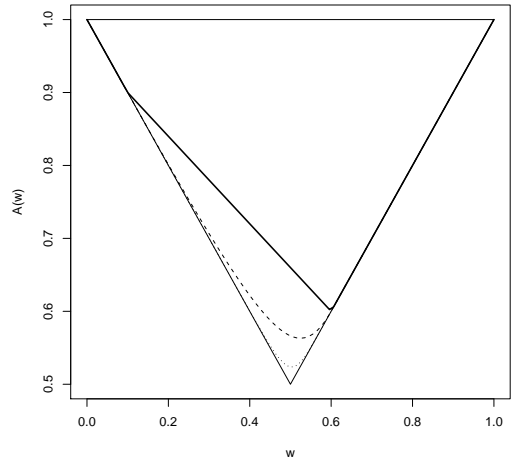

FIG 4. *We fix $a = 0.10$ and $b = 0.60$ but with varieties of $s = 1$ (solid line), $s = 2$ (dashed line) and $s = 5$ (dotted line) values.*

## 7. Estimating $c$ Using a Non-parametric Estimator of Dependence Function

In this section we show how to estimate $c$ using Pickands' non-parametric estimator of the dependence function. We show also that the estimated $c$ value obtained from Pickands' estimator is always greater than the true value of $c$. As $X < Y$ we know that

$$A(w) = 1 - w \quad \text{for} \quad 0 \leq w \leq c \tag{7.1}$$

for some $0 \leq c \leq 0.5$. We wish to estimate $c$ using Pickands' estimator defined in equation (2.5), the modified Pickands estimator by (5)

$$\widehat{A}_p(w) = \frac{n}{\sum_{i=1}^n \min\left(\frac{X}{w\overline{X}}, \frac{Y}{(1-w)\overline{Y}}\right)},$$

For $\widehat{A}_p(w) = 1 - w$ for given $0 < w < 0.5$ then $\frac{X_i}{w\overline{X}} > \frac{Y_i}{(1-w)\overline{Y}}$ $\forall i, i = 1, \ldots, n$. Let $W_i = \frac{X_i}{X_i + Y_i}$, then

$$\frac{X_i}{w\overline{X}} > \frac{Y_i}{(1-w)\overline{Y}} \Rightarrow \frac{W_i}{w\overline{X}} > \frac{1 - W_i}{(1-w)\overline{Y}}, \quad \forall i, i = 1, \ldots, n.$$

$$\Rightarrow W_i > \frac{w\overline{X}}{(1-w)\overline{Y} + w\overline{X}}$$

$$= \frac{1}{1 + \frac{\overline{Y}}{\overline{X}}\left(\frac{1}{w} - 1\right)} \quad \forall i, i = 1, \ldots, n.$$

As $n \to \infty$ $\overline{X} \to 1$ and $\overline{Y} \to 1$ so $\frac{1}{1+\frac{\overline{Y}}{\overline{X}}(\frac{1}{w}-1)} \longrightarrow w$. This leads to our conclusion that

$$\Rightarrow W_i > w \quad \forall i, i = 1, \ldots, n.$$

So if $A(w)$ satisfies equation (7.1), modified estimator satisfies

$$c < \widehat{c}_{\text{Pickands}} = \min(W_1, \ldots, W_n).$$

We recognise this is biased as $c < \widehat{c}_{Pickands}$ always.

## 8. Simulation study

We simulate 50 repetitions from $(Z_x, Z_y) \sim BEVM(s, \mu_x(t), \sigma_x, \mu_y(t), \sigma_y, \xi)$ where $s = 2$, $\mu_x(t) = 100 - 40t$, $\mu_y(t) = 150 - 40t$, $\sigma_x = 4$, $\sigma_y = 2$ and $\xi = 0.2$. This distribution satisfies Condition 1 of Theorem 2 in Chapter 6 with the marginal distribution for $Z_x$ is $\text{GEVM}(\mu_x(t), \sigma_x, \xi)$ and $Z_y$ is $\text{GEVM}(\mu_y(t), \sigma_y, \xi)$. These data are non-stationary with a trend in mean i.e. $\mu_x(t) = 100 - 40t$ and $\mu_y(t) = 150 - 40t$ which later are estimated using penalized log-likelihood, i.e. $\hat{g}_x$ and $\hat{g}_y$ respectively.

Figure 5 (left) shows realisations of $Z_x$ and $Z_y$ against time showing clearly that $Z_x < Z_y$ for all realisations. Figure 5 (right) shows that there is a dependence relationship for $Z_x$ and $Z_y$ after de-trending (removing the trend to give IID residuals).

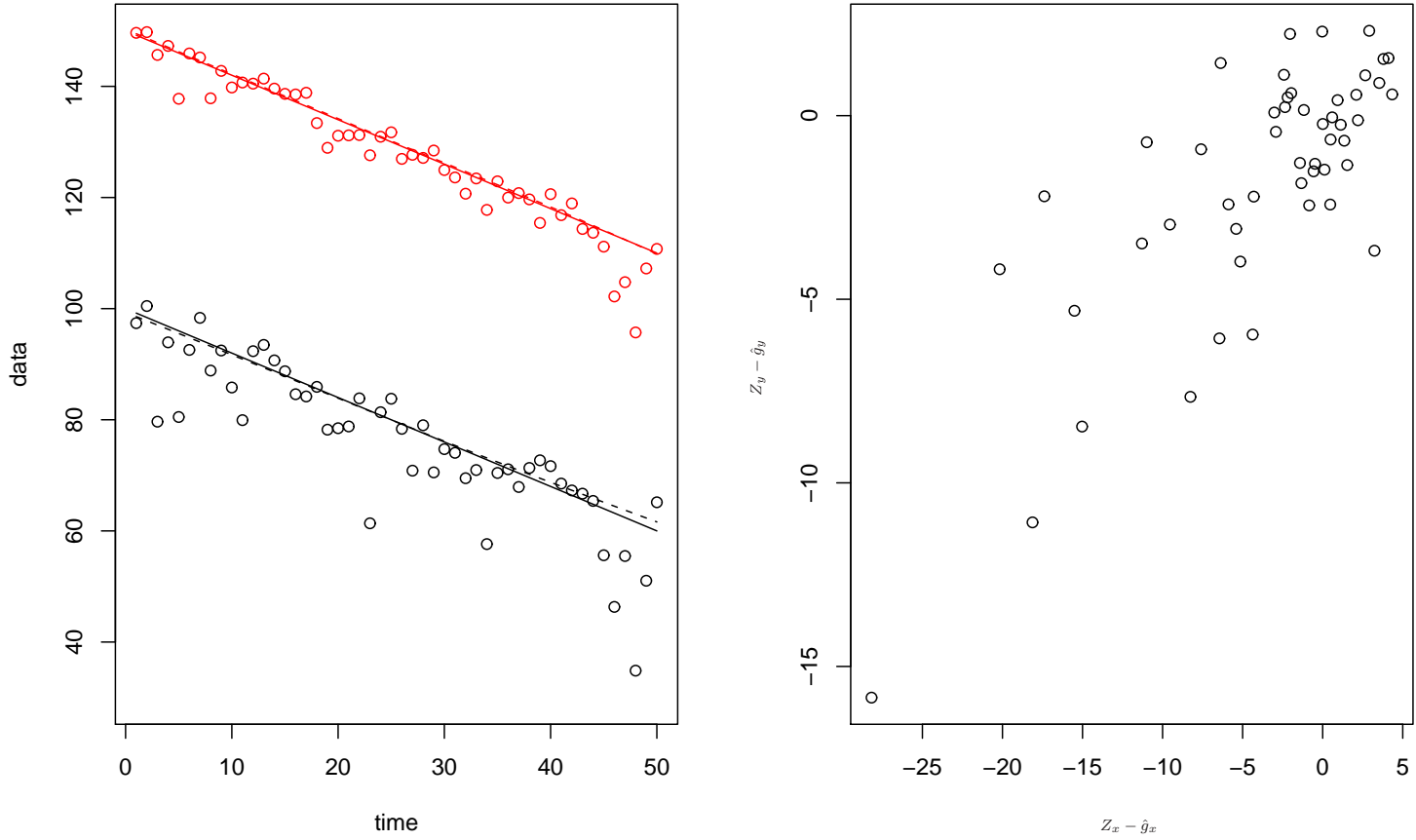

Fig 5. *Left plot is a scatter-plot for marginal data of $Z_x$ (black circle) with trend $\mu_x(t) = 100 - 40t$ (black line) and $\hat{g}_x$(black dashed line); and $Z_y$ (red circles) with trend $\mu_y(t) = 150 - 40t$ (red line) and $\hat{g}_y$(red dashed line). Right plot is $Z_x - \hat{g}_x$ versus $Z_y - \hat{g}_y$.*

TABLE 1
*The convergence values for parameters s, $\mu_x$, $\mu_y$ and $\xi$ from the simulation data.*

| s | $\sigma_x$ | $\sigma_y$ | $\xi$ |
|---|---|---|---|
| 1.68 | 3.71 | 1.85 | 0.214 |

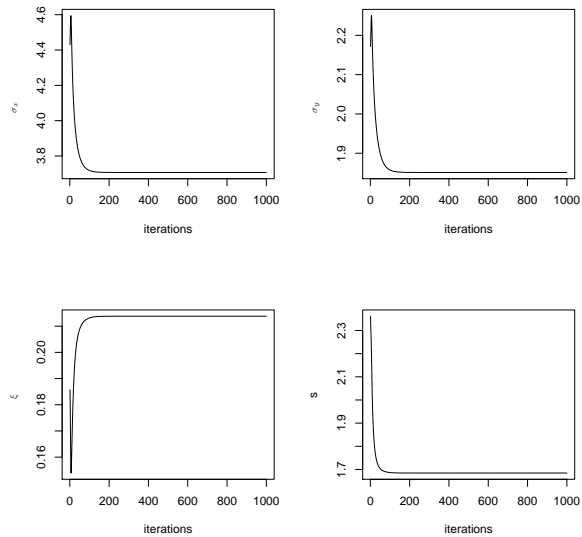

FIG 6. *The plot for converged parameter values for $\sigma_x$, $\sigma_y$, $\xi$ and s using penalised log-likelihood with smoothing parameter $\lambda_x = 1000$ and $\lambda_y = 1000$ for 1000 iterations.*

From Figure 6 we see the convergence for iterative estimates of the parameters $s$, $\sigma_x$, $\sigma_x$ and $\xi$, with final estimates tabulated in Table 1. We have the $\hat{g}_x$ and $\hat{g}_y$ represented by dashed lines in Figure 5 (left).

Figure 7 shows that there is strong dependence between $Z_x$ and $Z_y$ as all dependence functions are closer to the bounding curve $\max(w, 1-w)$ for $w \in [\hat{c}, 1]$ with values for each $\hat{c}$ for dependence functions.

From the working of the theory, the density measure $h_{RL}$ is defined on $\hat{c} = 0.030 < \hat{c}_{Pickands} = 0.078$ for true marginal parameter values and $\hat{c} = 0.052 < \hat{c}_{Pickands} = 0.070$ for estimated marginal parameter values from the simulation data.

Before assessing the model fit using P-P and Q-Q plots we need to transform the data to Exponential margins using equation (7.37) discussed in Subsection 7.5.2. By removing the marginal aspect in the data, the P-P and Q-Q plots in Figure 8 show that the data give a good fit to the BEVM of the restricted logistic function for the simulation data in Figure 5 when there is assumed to be dependence.

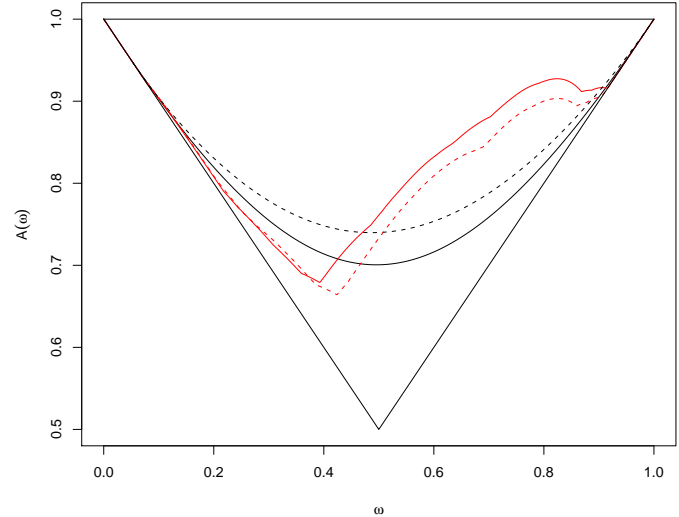

Fig 7. *Simulation data: Pickands dependence function estimates, $A(w)$ based on annual minima of Figure 5. Non-parametric estimates by Pickands using true parameter values (red line) and using estimated parameter values (dashed red line). Parametric estimates for restricted logistic function using true parameter values (black line) and using estimated parameter values (dashed black line).*

## 9. Discussion

When $X$ and $Y$ are ordered with the special constraint of stochastic ordering $X < Y$, it imposes som elimitations to the model of ordered bivariate extreme values. The logistics function defined by (6, Nadarajah) and (10, Tawn) was used to obtain a modify function called a restricted logistic function to model the dependence between variables. We proposed also some extensions of the dependence function of Nadarajah more generally. One of the important results is that when ordering exist, i.e. $X < Y$, we found that the density measure $h(w) = 0$, for $w \in [0, c = \max_{x \in R_Y} C(x,x)]$. This led to some other result, i.e. the restricted logistic dependence function (a modification version of asymetrical logistic function), restricted exponential function, $V_{RL}(x,y)$.

We showed through a simulation study that the ordered bivariate moel fits well when the constraint $X < Y$ is imposed. The simulation study with selection of parallel linear trends for $X$ and $Y$ showed that this approach can be applied to real data.

## Acknowledgment

I would like to thanks Professor Jonathan Tawn for supervisoring me on this research when I read my PhD at Lancaster University, United Kingdom from 2003 to 2007.

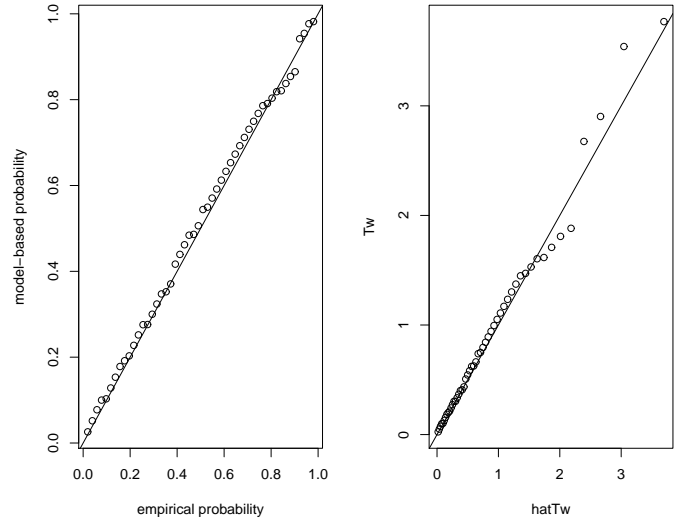

FIG 8. *The P-P (left) and Q-Q (right) plots for the simulation data of $(Z_x, Z_y) \sim BEVM(s = 2, \mu_x(t) = 100 - 40t, \sigma_x = 4, \mu_y(t) = 150 - 40t, \sigma_y = 2, \xi = 0.2)$ which followed to Condition 1 of Theorem 2 in Chapter 6 with 50 repetitions.*

**Proof** of equation (5.4).

The measure distribution function, $H_a(w)$, is derived as follows

$$H_a(w) = \int_0^w h_a(t)dt$$
$$= 1 - 2c + [w^{s-1} - (1-2c)^s(1-w)^{s-1}][w^s + (1-2c)^s(1-w)^s]^{\frac{1}{s}-1}.$$

We get the integral of $H_a(t)$ for $0 < t < w$ and $0 < t < 1$ as follows

$$\int_0^w H_a(t)dt = -(1-2c)(1-w) + [w^s + (1-2c)^s(1-w)^s]^{\frac{1}{s}},$$
$$\int_0^1 H_a(t)dt = 1,$$

and $H_a(0) = 0$ and $H_a(1) = 2(1-c)$.

The exponential measure function for restricted logistic function we simplify as follows

$$V_{RL}(x,y) = y\int_0^{D_+}(1-w)h_a(w)dw + \frac{cx}{1-c}\int_{D_+}^1 h_a(w)dw + x\int_{D_+}^1 wh_a(w)dw,$$

where $D_+ = \max(D, 0)$ with $0 \le D_+ \le 1$ and

$$D = \frac{y}{x+y} - \left(\frac{c}{1-c}\right)\frac{x}{x+y} = 1 - \left(\frac{1}{1-c}\right)\frac{x}{x+y} = \frac{1}{1-c}\left(\frac{y}{x+y} - c\right).$$

When $c < \frac{y}{x+y} < 1$ then $D_+ = D$, the derivation of $V_{RL}$ is using integration by parts, letting $u = 1 - w$ and $dv = h_a(w)dw \Rightarrow du = -dw$ and $v = H_a(w)$, we get

$$\int_0^{D_+} (1-w)h_a(w)dw = L(D_+) + M(D_+)P(D_+)(1-D_+) \qquad (9.1)$$

where

$$L(D_+) = \left[D_+^s + (1-2c)^s(1-D_+)^s\right]^{\frac{1}{s}}, \qquad (9.2)$$

$$M(D_+) = \left[D_+^s + (1-2c)^s(1-D_+)^s\right]^{\frac{1}{s}-1}, \qquad (9.3)$$

$$P(D_+) = D_+^{s-1} - (1-2c)^s(1-D_+)^{s-1}. \qquad (9.4)$$

Letting $u = w$ and $dv = h_a(w)dw \Rightarrow du = dw$ and $v = H_a(w)$,

$$\int_{D_+}^1 w h_a(w)dw = L(D_+) - D_+ M(D_+)P(D_+) \qquad (9.5)$$

$$\int_{D_+}^1 h_a(w)dw = 1 - M(D_+)P(D_+). \qquad (9.6)$$

Using equations (9.1), (9.2), (9.4),(9.3), (9.5) and (9.6) then

$$V_{RL}(x,y) = (y+x)\,L(D_+) + \frac{cx}{(1-c)},$$

as

$$yM(D_+)P(D_+) - \frac{cx}{(1-c)}M(D_+)P(D_+) = yD_+M(D_+)P(D_+) + xD_+M(D_+)P(D_+).$$

When $y = 0$ then $D_+ = 0$.

Then the exponential measure function is

$$V_{RL}(x,0) = x.$$

When $x = 0$ then $D_+ = 1$. Then the exponential measure function is $V_{RL}(0,y) = y$.

When $c = 0$, then $D_+ = \frac{y}{x+y}$ we get $V_{RL}(x,y) = (y^s + x^s)^{\frac{1}{s}}$. When $0 \le \frac{y}{x+y} \le c$ then $D_+ = 0$ we get

$$V_{RL}(x,y) = 2cx + (1-2c)x = x. \qquad (9.7)$$

Checking for $c = 0$ then when $x = 0$, we get $V(0,y) = y$. When $y = 0$, we get $V(x,0) = x$. Another way of writing the exponential measure function is by including the dependence function, $A$,

$$V_{RL}(x,y) = (x+y)\,A_{RL}\left(\frac{x}{x+y}\right). \qquad (9.8)$$

From equation (9.7), for $0 < w = \frac{y}{x+y} < c$

$$x = (x+y) A_{RL}\left(\frac{x}{x+y}\right)$$

$$\Rightarrow A_{RL}(w) = \frac{x}{x+y} = 1 - w.$$

When $c < w = \frac{y}{x+y} < 1 \Rightarrow D_+ = D$ from equation (9.2) and relation with equation (9.8) we obtain

$$A_{RL}(w) = \frac{1}{1-c}\left\{c(1-w) + [(w-c)^s + (1-2c)^s(1-w)^s]^{\frac{1}{s}}\right\}$$

Checking when $w = c$ then

$$A_{RL}(c) = \frac{1}{1-c}[c(1-c) + (1-2c)(1-c)] = c + 1 - 2c = 1 - c.$$

Considering all possible values for $D_+$ the generalized form of the exponential measure for restricted logistic is

$$V_{RL}(x,y) = \begin{cases} x & \text{for } 0 < \frac{y}{x+y} \leq c \\ \frac{1}{1-c}\left(\{[(1-c)y - cx]^s + (1-2c)^s x^s\}^{\frac{1}{s}} + cx\right) & \text{for } c < \frac{y}{x+y} < 1. \end{cases}$$

Using equation (9.8) we can easily get the restricted dependence function of $A_{RL}(x,y)$ where for $0 < \omega \leq c$

$$A_{RL}(x,y) = 1 - \omega$$

and for $c < \omega < 1$, we get

$$A_{RL}(x,y) = \frac{1}{1-c}\left\{c(1-\omega) + \left[(\omega-c)^{\frac{1}{\alpha}} + (1-2c)^{\frac{1}{\alpha}} + (1-2c)^{\frac{1}{\alpha}}(1-\omega)^{\frac{1}{\alpha}}\right]^\alpha\right\}.$$

From (3) the relationship between an exponential measure, $V$, and a density measure function, $h$, involving two variables of Fréchet margins denoted by "F" where $x = \frac{1}{x_F}$ and $y = \frac{1}{x_F}$ is

$$\frac{\partial^2 V_{RL}(x_F, y_F)}{\partial x_F \partial y_F} = -\frac{1}{(x_F + y_F)^3} h_{RL}\left(\frac{x_F}{x_F + y_F}\right).$$

By changing the equation (5.4) to Fréchet margins we get $V_{RL}(x_F, y_F)$ is equal to

$$V_{RL} = \begin{cases} \frac{1}{x_F} & \text{for } 0 < \frac{x_F}{x_F + y_F} \leq c \\ \frac{1}{1-c}\left\{\left[\left(\frac{1-c}{y_F} - \frac{c}{x_F}\right)^s + (1-2c)^s\left(\frac{1}{x_F}\right)^s\right]^{\frac{1}{s}} + \frac{c}{x_F}\right\} & \text{for } c < \frac{x_F}{x_F + y_F} < 1. \end{cases}$$

(9.9)

From equation (9.9), $\frac{\partial^2 V_{RL}(x_F, y_F)}{\partial x_F \partial y_F} = 0$. As $y_F \to \infty$, $\frac{x_F}{x_F+y_F} \to 0$ we know $h_{RL}(0) = 0$. We get for $0 < \frac{x_F}{x_F+y_F} \leq c$,

$$\frac{\partial^2 V_{RL}(x_F, y_F)}{\partial x_F \partial y_F} = -\frac{1}{(x_F+y_F)^3} h_{RL}\left(\frac{y_F}{x_F+y_F}\right) = 0.$$

From equation (9.9) checking for $c < \frac{x_F}{x_F+y_F} < 1$, we get

$$\begin{aligned}\frac{\partial^2 V_{RL}(x_F, y_F)}{\partial x_F \partial y_F} &= -\frac{1}{(y_F x_F)^3}(1-c)(s-1)(1-2c)^s \left\{\frac{1}{x_F}\left(\frac{1-c}{y_F} - \frac{c}{x_F}\right)\right\}^{s-2} \\ &\quad \times \left[\left(\frac{1-c}{y_F} - \frac{c}{x_F}\right)^s + (1-2c)^s\left(\frac{1}{x_F}\right)^s\right]^{\frac{1}{s}-2} \quad (9.10)\\ &= -(s-1)(1-c)(1-2c)^s \left\{y_F\left(x_F(1-c) - cy_F\right)\right\}^{s-2} \\ &\quad \times \left[(x_F(1-c) - cy_F)^s + (1-2c)^s y_F^s\right]^{\frac{1}{s}-2}.\end{aligned}$$

$$\begin{aligned}h_{RL}\left(\frac{x_F}{x_F+y_F}\right) &= (s-1)(1-c)(1-2c)^s \left\{\left(\frac{x_F}{x_F+y_F} - c\right)\left(\frac{y_F}{x_F+y_F}\right)\right\}^{s-2} \\ &\quad \times \left\{\left(\frac{x_F}{x_F+y_F} - c\right)^s + (1-2c)^s\left(\frac{y_F}{x_F+y_F}\right)^s\right\}^{\frac{1}{s}-2} \\ &= \left(\frac{1}{x_F+y_F}\right)^{-3}(s-1)(1-c)(1-2c)^s \left\{y_F\left[x_F(1-c) - cy_F\right]\right\}^{s-2} \\ &\quad \times \left\{[x_F(1-c) - cy_F]^s + (1-2c)^s y_F^s\right\}^{\frac{1}{s}-2} \\ &= -(x_F+y_F)^3 \frac{d^2 V_{RL}(x_F, y_F)}{dxdy}\end{aligned}$$

$$\Rightarrow \frac{d^2 V_{RL}(x_F, y_F)}{dx_F dy_F} = -\frac{1}{(x_F+y_F)^3} h_{RL}\left(\frac{x_F}{x_F+y_F}\right).$$